\documentclass[a4paper,11pt]{article}
\usepackage{latexsym}
\usepackage{amsfonts,amssymb}
\usepackage{amscd}
\usepackage[dvips]{graphicx}
\usepackage{color}
\font\erm = cmr9

\author{M.~Dziemia\'nczuk}
\title{On Cobweb Posets and Discrete F-Boxes Tilings}

\newtheorem{defn}{Definition}
\newtheorem{problem}{Problem}
\newtheorem{theoremn}{Theorem}
\newtheorem{observ}{Observation}
\newtheorem{proposition}{Proposition}
\newtheorem{notation}{Notation}
\newtheorem{corol}{Corollary}

\newtheorem{fact}{Fact}

\newcommand{\layer}[2]{\nobreak{\langle\Phi_{#1}\!\to\!\Phi_{#2}\rangle}}

\newcommand{\fnomial}[2]{ {{#1} \choose {#2}}_{\!\!F} }
\newcommand{\fnomialF}[3]{ {{#1} \choose {#2}}_{\!\!#3} }

\newcommand{\Nat}{\mathbb{N}}
\newcommand{\NatZero}{\mathbb{N}\cup\{0\}}

\begin{document}

\begin{center}
{\textsc {\Large On Cobweb Posets\\ and Discrete F-Boxes Tilings}}  \\ 

\vspace{0.5cm} 
Maciej Dziemia\'nczuk

\vspace{0.5cm}
{\erm
	Institute of Informatics, University of Gda\'nsk \\
	PL-80-952 Gda\'nsk, Wita Stwosza 57, Poland\\
	e-mail: mdziemianczuk@gmail.com\\
}

\end{center}

\begin{abstract}

$F$-boxes defined in \cite{akkmd2} as hyper-boxes in $N^{\infty}$ discrete space were applied here for the geometric description of the cobweb posetes Hasse diagrams tilings. The $F$-boxes edges sizes are taken to be values of terms of natural numbers' valued sequence $F$.  The problem of partitions of hyper-boxes represented by graphs into blocks of special form is considered and these are to be called $F$-tilings.

The proof of such tilings' existence for certain sub-family of admissible sequences $F$ is delivered. The family of $F$-tilings which we consider here includes among others  $F$ = Natural numbers, Fibonacci numbers, Gaussian integers with their corresponding $F$-nomial (Binomial, Fibonomial, Gaussian) coefficients  as it is persistent typical for combinatorial interpretation of such tilings originated from Kwa\'sniewski cobweb posets tiling problem . 

Extension of this tiling problem onto the general case multi $F$-nomial coefficients is here proposed. Reformulation of the present cobweb tiling problem into a clique problem of a graph specially invented for that purpose - is proposed here too. To this end we illustrate the area of our reconnaissance by means of the Venn type map of various cobweb sequences families.

\vspace{0.4cm}
\noindent AMS Classification Numbers: 05A10, 05A19, 11B83, 11B65

\vspace{0.2cm}
\vspace{0.3cm}

\noindent Affiliated to The Internet Gian-Carlo Polish Seminar: \\
\noindent \emph{http://ii.uwb.edu.pl/akk/sem/sem\_rota.htm}, \\
Article \textbf{No7}, April  2009, 15  April  2009, \\
(302 anniversary of Leonard Euler's birth)

\end{abstract}

\vspace{0.3cm}
\section{Introduction}

The \emph{Kwa\'sniewski upside-down} notation from \cite{akk4} (see also \cite{akk1,akk2}) is being here taken for granted. For example $n$-th element of sequence $F$ is $F_n \equiv n_F$, consequently $n_F! = n_F\cdot(n-1)_F\cdot...\cdot 1_F$ and a set $[n_F] = \{1,2,...,n_F\}$ however $[n]_F=\{1_F,2_F,...,n_F\}$. More about effectiveness of this notation see references in \cite{akk4} and  Appendix ``\emph{On upside-down notation}'' in \cite{akkmd2}.

Throughout this paper we shall consequently use $F$ letter for a sequence of positive integers i.e. $F\equiv\{n_F\}_{n\geq 0}$ such that $n_F\in\Nat$ for any $n\in\NatZero$.

\subsection{Discrete $m$-dimensional $F$-Box}

Let us define discrete $m$-dimensional $F$-box with edges sizes designated by natural numbers' valued sequence $F$ as described below.
These $F$-boxes from \cite{akkmd2} where invented as a response to \emph{Kwa\'sniewski cobweb tiling} problem posed in \cite{akk1} (Problem 2 therein) and his question about visualization of this phenomenon.

\begin{defn}
Let $F$ be a natural numbers' valued sequence $\{n_F\}_{n\geq 0}$ and $m,n\in\Nat$ such that $n\geq m$. Then a set $V_{m,n}$ of points $v=(v_1,...,v_m)$ of discrete $m$-dimensional space $\Nat^m$ given as follows
\begin{equation}
	V_{m,n} = [k_F] \times [(k+1)_F] \times ... \times [n_F]
\end{equation}
where $k=n-m+1$ and $[s_F] = \{1,2,...,s_F\}$ is called $m$-dimensional $F$-box. 
\end{defn}

\begin{figure}[ht]
\begin{center}
	\includegraphics[width=100mm]{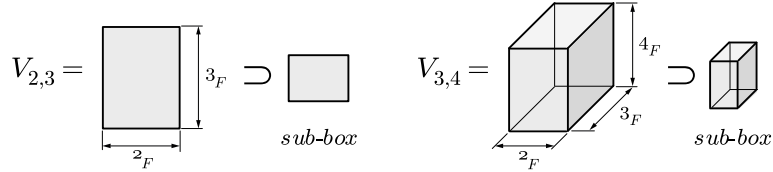}
	\caption{$F$-Boxes $V_{2,3}$ and $V_{3,4}$ with sub-boxes.}
\end{center}
\end{figure}

\begin{figure}[ht]
\begin{center}
	\includegraphics[width=110mm]{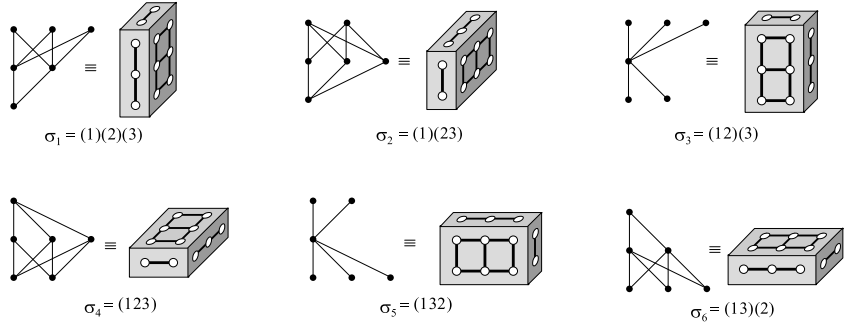}
	\caption{$F$-Boxes of the form $\sigma V_3$ and cobweb blocks $\sigma P_3$. \label{fig:Tiles3D}}
\end{center}
\end{figure}

In the case of $n = m$ we write for short $V_{m,m} \equiv V_{m}$. Assume that we have a $m$-dimensional box $V_{m,n} = W_1\times W_2\times...\times W_m$. Then a set $A = A_1\times A_2\times ... \times A_m$ such that
$$
	A_s \subset W_s, \qquad |A_s|>0, \qquad s=1,2,...,m;
$$
is called \emph{$m$-dimensional sub-box of $V_{m,n}$}. Moreover, if for $s=1,2,...,m$ these sets $A_s$ satisfy the following
$$
	|A_s| = (\sigma\cdot s)_F
$$
for any permutation $\sigma$ of set $\{1_F,2_F,...,m_F\}$ then $A$ is called \emph{$m$-dimensional sub-box of the form $\sigma V_m$}. Compare with Figure \ref{fig:Tiles3D}.

\vspace{0.2cm}
Note, that the permutation $\sigma$ might be understood here as an orientation of sub-box's position in the box $V_{m,n}$. Any two sub-boxes $A$ and $B$ are disjoint if its sets of points are disjoint i.e. $A\cap B = \emptyset$.

\vspace{0.2cm}
The number of points $v=(v_1,...,v_m)$ of $m$-dimensional box $V_{m,n}$ is called \emph{volume}.
It it easy to see that the \emph{volume} of $V_{m,n}$ is equal to

\begin{equation}
	|V_{m,n}| = n_F\cdot (n-1)_F \cdot ... \cdot (n-m+1)_F = n^{\underline{m}}_F
\end{equation}
while for  $m = n$ 
\begin{equation}
	|V_m| = |\sigma V_m| = m_F \cdot (m-1)_F \cdot ... \cdot 1_F =  m_F!
\end{equation}

\subsection{Partition of discrete $F$-boxes}

Let us consider $m$-dimensional $F$-box $V_{m,n}$. A finite collection of $\lambda$ pairwise disjoint sub-boxes $B_1,B_2,...,B_\lambda$ of the volume equal to $\kappa$ is called \emph{$\kappa$-partition} of $V_{m,n}$ if their set union of gives the  whole box $V_{m,n}$ i.e.
\begin{equation}
	\bigcup_{1\leq j \leq \lambda} B_j = V_{m,n}, \qquad |B_i| = \kappa,\qquad  i=1,2,...,\lambda.
\end{equation}

\vspace{0.2cm}

\noindent \textbf{Convention.} In the following, we shall deal  only with these $\kappa$-partition of $m$-dimensional boxes $V_{m,n}$, which volume $\kappa$ of sub-boxes is equal to the volume of box $V_m$ i.e. $\kappa = |V_m|$.

Of course the box $V_{m,n}$ has $\kappa$-partition  \textit{not for all} $F$ - sequences \cite{md1}. Therefore we introduce the name: \emph{$F$-admissible} sequence which means that $F$  satisfies the necessary and sufficient conditions for the box $V_{m,n}$ to have $\kappa$-partitions.  
In order to proceed let us recall first what follows.

\begin{defn}[\cite{akk1,akk2}]\label{def:fnomial}
Let $F$ be a natural numbers' valued sequence $F=\{n_F\}_{n\geq 0}$. Then $F$-nomial coefficient is identified with the symbol
\begin{equation}
	\fnomial{n}{m} = \frac{n_F!}{m_F!(n-m)_F!} = \frac{n_F^{\underline{m}}}{m_F!}
\end{equation}
where $n_F^{\underline{0}} = 0_F! = 1$.
\end{defn}

\begin{defn}[\cite{akk1,akk2}] \label{def:admissible}
A sequence $F$ is called admissible if, and only if for any $n,m\in \NatZero$ the value of $F$-nomial coefficient is natural number or zero i.e.
\begin{equation}\label{eq:admissible}
	\fnomial{n}{m} \in \NatZero
\end{equation} 
while $n\geq m$ else is zero.
\end{defn}

\vspace{0.2cm}
Recall now also a combinatorial interpretation of the $F$-nomial coefficients in $F$-box reformulated form (consult Remark 5 in \cite{akk4} and \cite{akkmd2}). And note: these coefficients encompass among others Binomial, Gaussian and Fibonomial coefficients.

\begin{fact}[Kwa\'sniewski \cite{akk1,akk2}]
Let $F$ be an admissible sequence. Take any $m,n\in\Nat$ such that $n\geq m$, then the value of $F$-nomial coefficient $\fnomial{n}{m}$ is equal to the number of sub-boxes that constitute a  $\kappa$-partition of $m$-dimensional $F$-box $V_{m,n}$ where $\kappa = |V_m|$.
\end{fact}

\noindent {\it{\textbf{Proof.}}}
This proof comes from Observation 3 in \cite{akk1,akk2} and was adopted here to the language of discrete boxes. Let us consider $m$-dimensional box $V_{m,n}$ with  $|V_{m,n}| = n^{\underline{m}}_F$. The volume of sub-boxes is equal to $\kappa = |V_m| = m_F!$. 
Therefore the number of sub-boxes is equal to
$$
	\frac{n^{\underline{m}}_F}{m_F!} = \fnomial{n}{m}
$$
From definition of $F$-admissible sequence we have that the above is natural number. Hence the thesis $\blacksquare$

\vspace{0.4cm}
While considering any $\kappa$-partition of certain $m$-dimensional box we only assume that sub-boxes \textbf{have the same volume}. In the next section we shall take into account these partitions which sub-boxes have additionally established structure.

\subsection{Tiling problem}

Now, special $\kappa$-partitions of discrete boxes are considered. 
Namely, we deal with only these partitions of $m$-dimensional box $V_{m,n}$ which all sub-boxes \textbf{are of the form} $V_m$.

\begin{defn}
Let $V_{m,n}$ be a $m$-dimensional $F$-box. Then any $\kappa$-partition into sub-boxes of the form $V_m$ is called tiling of $V_{m,n}$.
\end{defn}

It was shown in \cite{md1} that just the admissibility condition (\ref{eq:admissible}) is not sufficient for the existence a tiling for any given $m$-dimensional box $V_{m,n}$. Kwa\'sniewski in his papers \cite{akk1,akk2} posed the following problem called \emph{Cobweb Tiling Problem}, which was a starting point of the research with  results being reported in the presents note.

\begin{problem} [Tiling]
Suppose now that $F$ is an admissible sequence. Under which conditions any $F$-box $V_{m,n}$ designated by sequence $F$ has a tiling? Find effective characterizations and/or find an algorithm to produce these tilings.
\end{problem}

\begin{figure}[ht]
\begin{center}
	\includegraphics[width=80mm]{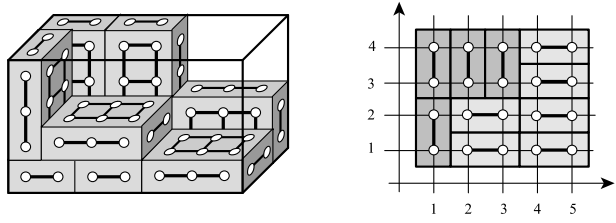}
	\caption{Sample 3D and 2D tilings.}
\end{center}
\end{figure}

In the next sections we propose certain family $\mathcal{T}_\lambda$ of sequences $F$. Then we prove that any $F$-box $V_{m,n}$, where $m,n\in\Nat$ designated by $F\in\mathcal{T}_\lambda$ has a tiling with giving a construction of it.

\subsection{Cobweb representation}

In this section we recall \cite{akkmd2} that discrete $F$-boxes $V_{m,n}$ are unique codings  representing \emph{Cobwebs}, introduced by Kwa\'sniewski \cite{akk1,akk2} as a special graded posets. Any poset might be represented as a Hasse digraph and this approach to tiling problem will be used throughout the paper.

Next we shall consider partitions of $m$-dimensional boxes as a partitions of cobwebs with $m$ levels into sub-cobwebs called blocks. In the following we quote some necessary notation of \emph{Cobwebs} adopted to the tiling problem. For more on \emph{Cobwebs} see source papers \cite{akk1,akk2,akk4} and references therein.

\begin{defn} 
Let $F$ be a natural numbers' valued sequence. Then a simple graph $\langle V, E\rangle$, such that $V = \bigcup_{k \leq s \leq n} \Phi_s$ and
\begin{equation}
	 E=\Big\{ \{ u,v \} : u\in\Phi_s \wedge v\in\Phi_{s+1} \wedge k\geq s < n \Big\}
\end{equation}
where $\Phi_s = \{1,2,...,s_F\}$ is called cobweb layer $\layer{k}{n}$.
\end{defn}

\begin{figure}[ht]
\begin{center}
	\includegraphics[width=50mm]{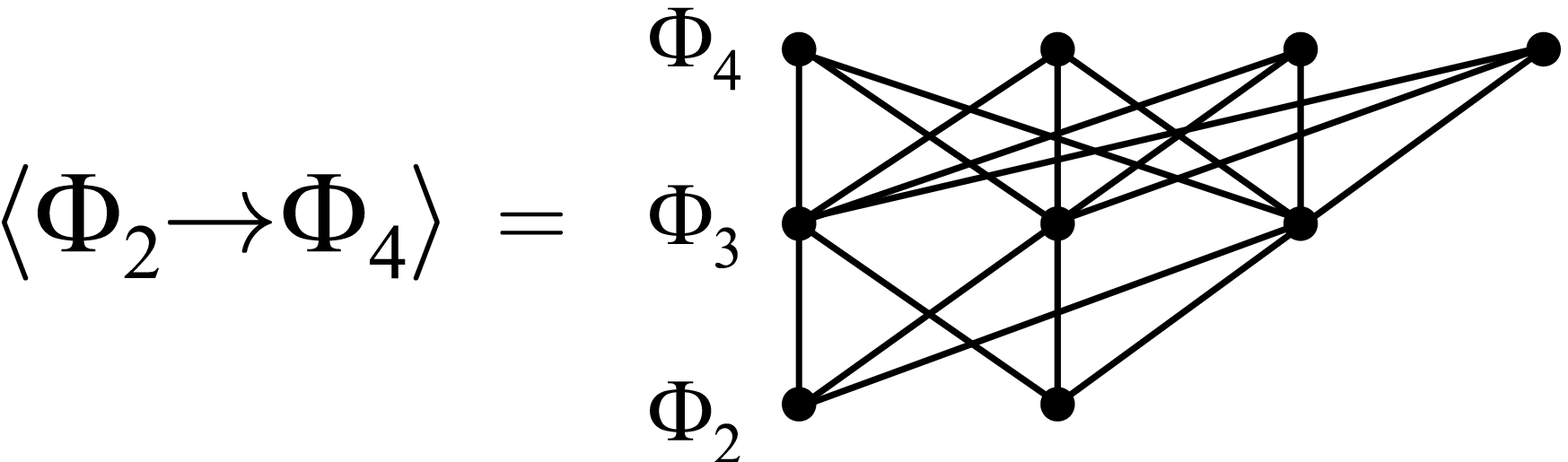}
	\caption{Cobweb layer $\layer{2}{4}$ designated by $F$=Natural numbers. \label{fig:Layer}}
\end{center}
\end{figure}

Suppose that we have a cobweb layer $\layer{k}{n}$ of $m$ levels $\Phi_s$, where $m = n-k+1$. Then any cobweb layer $\langle\phi_1\to\phi_m\rangle$ of $m$ levels $\phi_s$ such that
\begin{equation}
	\phi_s \subseteq \Phi_s, \qquad |\phi_s| = s_F,  \qquad s=1,2,...,m;
\end{equation}
is called \emph{cobweb block} $P_m$ of layer $\layer{k}{n}$.

\begin{figure}[ht]
\begin{center}
	\includegraphics[width=80mm]{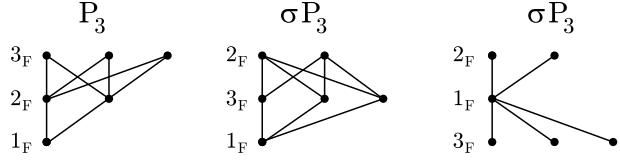}
	\caption{Example of cobweb blocks $P_3$ and $\sigma P_3$. \label{fig:Blocks}}
\end{center}
\end{figure}

\vspace{0.2cm}
Additionally, one considers cobweb blocks obtained via permutation $\sigma$ of theirs levels' order as follows (Compare with Figure \ref{fig:Blocks}).

\begin{defn}
Let a cobweb layer $\layer{k}{n}$ with  $m$ levels $\Phi_s$ be given, where $m=n-k+1$. Then a cobweb block $P_m$ with $m$ levels $\phi_s$ such that
\begin{equation}
	\phi_s \subseteq \Phi_s, \qquad |\phi_s| = (\sigma\cdot s)_F, \qquad s=1,2,...,m;
\end{equation}
where $\sigma$ is a permutation of the set $\{1_F,2_F,...,m_F\}$ is called cobweb block of the form $\sigma P_m$.
\end{defn}

\begin{figure}[ht]
\begin{center}
	\includegraphics[width=50mm]{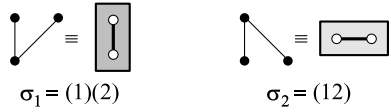}
	\caption{$F$-Boxes of the form $\sigma V_2$ and cobweb blocks $\sigma P_2$. \label{fig:Tiles2D}}
\end{center}
\end{figure}

While saying \emph{``a block $\sigma P_m$ of layer $\layer{k}{n}$''} we mean that the number of levels in block and layer is the same i.e. $m = n - k + 1$ and each of levels of block are non-empty subsets of corresponding levels in the layer.

Assume that we have a cobweb layer $\layer{k}{n}$. A path $\pi$ from any vertex at first level $\Phi_k$ to any vertex at the last level $\Phi_n$, such that
$$
	\pi = \{v_k, v_{k+1}, ..., v_{n}\}, \qquad v_s \in \Phi_s, \qquad s=k,k+1,...,n;
$$
is noted as a \emph{maximal-path $\pi$} of $\layer{k}{n}$.
In the same way we nominate  \emph{maximal-path} of cobweb block $\sigma P_m$.

\vspace{0.2cm}
Let $C_{max}(A)$ denotes a set of maximal-paths $\pi$ of cobweb block $A$. (Compare with [4]). Two cobweb blocks $A, B$ of layer $\layer{k}{n}$ are max-disjoint or disjoint for short (\cite{akk1,akk2}) if, and only if its sets of maximal-paths are disjoint i.e. $C_{max}(A) \cap C_{max}(B) = \emptyset$. 
The cardinality of set $C_{max}(A)$ is called \emph{size} of block $A$.

\begin{observ}[\cite{akkmd2}]
Let $F$ be a natural numbers' valued sequence and $k,n\in\Nat$. 
Then any $F$-box $V_{m,n}$ is uniquely represented by cobweb layer $\layer{k}{n}$ and vice versa i.e., 
\begin{equation}
	V_{m,n} \Leftrightarrow \layer{k}{n}.
\end{equation}
where $k = n-m+1$.
\end{observ}

\noindent {\it{\textbf{Proof.}}} Consider a cobweb layer $\layer{k}{n}$ of $m$ levels $\Phi$ and $m$-dimensional box $V_{k,n}$. Observe that any maximal-path $\pi = (v_1,v_2,...,v_m)$ of the layer corresponds to only one point $x = (x_1,x_2,...,x_m)$ of $m$-dimensional box $V_{m,n}$, and vice versa, i.e.
$$
	[s_F] \ni x_s \Leftrightarrow v_s \in [s_F],\qquad s=1,2,...,m;
$$
And the number of these maximal-paths and points is the same (Compare with \cite{akk4} and \cite{akkmd2}) i.e. 
$$
	|C_{max}(\layer{k}{n})| = |V_{m,n}|
$$
where $m=n-k+1$. $\blacksquare$

\begin{figure}[ht]
\begin{center}
	\includegraphics[width=100mm]{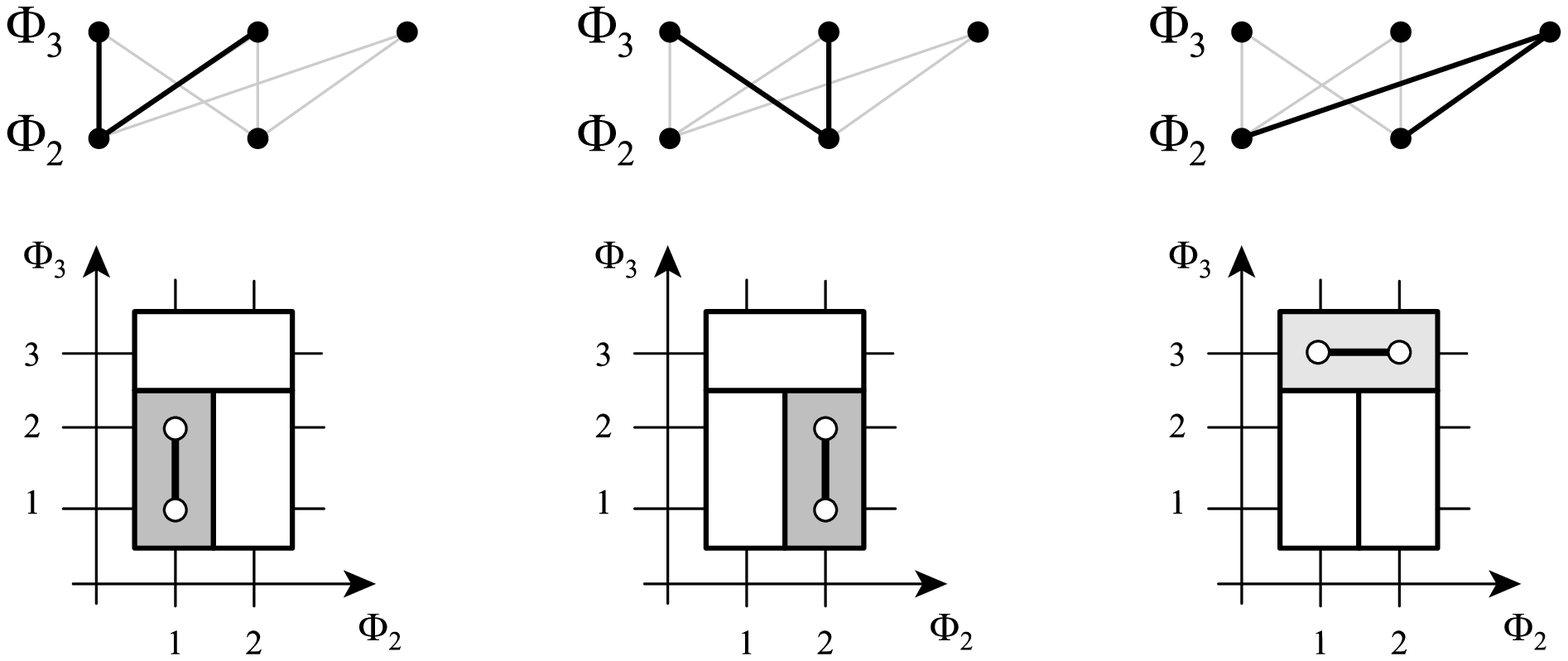}
	\caption{Correspondence between tiling of $F$-box $V_{3,4}$ and $\layer{3}{4}$.}
\end{center}
\end{figure}

\vspace{0.2cm}
Next, we draw terminology of $F$-boxes' partitions back to cobweb's language, used in the next part of this note.

Take any cobweb layer $\layer{k}{n}$ with $m$ levels. Then a set of $\lambda$ pairwise disjoint cobweb blocks $A_1,A_2,...,A_\lambda$ of $m$ levels such that its size is equal to $\kappa$ and the union of $C_{max}(A_1), C_{max}(A_2),...,C_{max}(A_\lambda)$ is equal to the set $C_{max}(\layer{k}{n})$ is called \emph{cobweb $\kappa$-partition}.
Finally, a $\kappa$-partition of layer $\layer{k}{n}$ with $m$ levels into cobweb blocks of the form $\sigma P_m$ is called \emph{cobweb tiling}.

\vspace{0.2cm}
\noindent Let us sum it up with the following Table \ref{tab:equiv}.

\renewcommand{\arraystretch}{1.3}

\begin{table}[ht]
\caption{Equivalent notation and terminology. \label{tab:equiv}}
\begin{center}
\begin{tabular}{|l |l | l |}
\hline
& \textbf{Cobwebs} & \textbf{$F$-boxes} \\
\hline 
1. & Maximal-path $(v_1,...,v_m) \in \layer{k}{n}$  &   Point $(x_1,...,x_m) \in V_{m,n}$ \\
2. & Cobweb layer $\layer{k}{n}$  &  $F$-box $V_{m,n}$ \\
3. & Cobweb block $\sigma P_m \subset \layer{k}{n}$  &  Sub-box $\sigma V_m \subset V_{m,n}$ \\
4. & Tiling of cobweb layer   &  Tiling of $F$-box \\
&  where k = n-m+1. &  \\
\hline
\end{tabular}
\end{center}
\end{table}


\section{Cobweb tiling sequences}

\begin{figure}[ht]
\begin{center}
	\includegraphics[width=100mm]{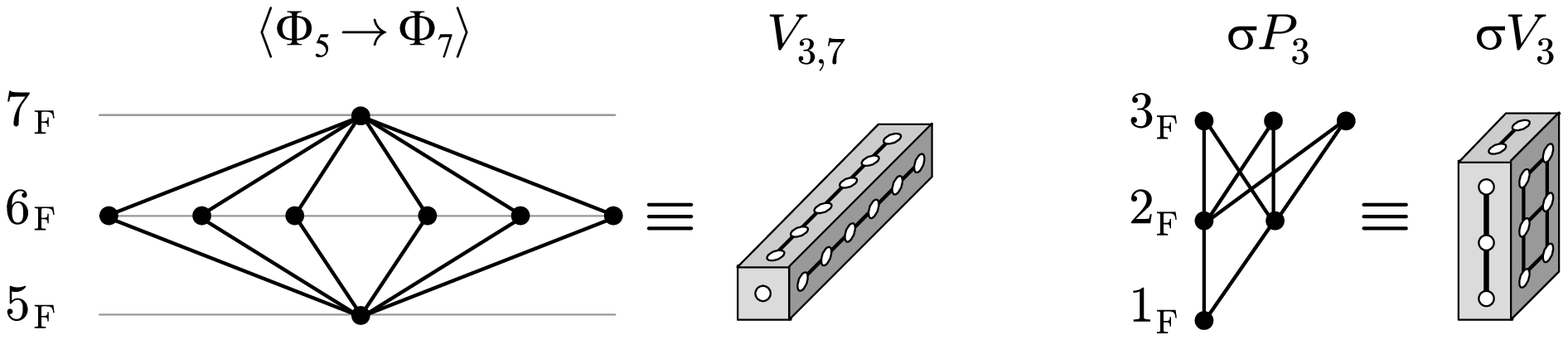}
	\caption{Layer $\layer{5}{7}$ that does not have tiling with blocks $\sigma P_3$. \label{fig:contr}}
\end{center}
\end{figure}

Recall that for some \emph{$F$-admissible} sequences there is no method to tile certain $F$-boxes $V_{m,n}$ or accordingly cobweb layers $\layer{k}{n}$ (no tiling property). For example see Figure \ref{fig:contr} that comes from \cite{md1}.
In the next part of this note, we define and consider \textbf{only} sequences \textbf{with tiling property}.

\begin{defn}
A cobweb admissible sequence $F$ such that for any $m,n\in\Nat$ the  cobweb layer $\layer{k}{n}$ has a tiling is called cobweb tiling sequence.
\end{defn}

Let $\mathcal{T}$ denotes the family of all cobweb tiling sequences. Characterization of whole family  $\mathcal{T}$ is still open problem. Nevertheless we define certain subfamily $\mathcal{T}_\lambda \subset \mathcal{T}$ of non-trivial cobweb tiling sequences. This family contains among others Natural and Fibonacci numbers, Gaussian integers and others.

\begin{notation}
Let $\mathcal{T}_\lambda$ denotes the family of natural number's valued sequences $F\equiv\{n_F\}_{n\geq 1}$ such that for any $n$-th term of $F$ satisfies the following holds
\begin{equation}\label{eq:TLambda}
	\forall~m,k\in\Nat,  \quad n_F = (m+k)_F = \lambda_K \cdot k_F \ + \ \lambda_M \cdot m_F  
\end{equation}
while $1_F\in\Nat$ and for certain coefficients $\lambda_K\equiv\lambda_K(k,m)\in\NatZero$ and $\lambda_M\equiv\lambda_M(k,m)\in\NatZero$.
\end{notation}

Note, coefficients $\lambda_K$ and $\lambda_M$ might be considered as a natural numbers' with zero valued infinite matrixes $\lambda_K \equiv [k_{ij}]_{i,j\geq 1}$ and $\lambda_M \equiv [m_{ij}]_{i,j\geq 1}$. Moreover the sequence $F\equiv\{n_F\}_{n\geq 0}$ is uniquely designated by these matrixes $\lambda_K, \lambda_M$ and first element $1_F\in\mathbb{N}$.

\begin{corol}\label{cor:lambdamult}
Let a sequence $F\in\mathcal{T}_\lambda$ with its coefficients' matrixes $\lambda_K,\lambda_M$ 
and a composition $\vec{\beta}=\langle b_1,b_2,...,b_k \rangle$ of number $n$ into $k$ nonzero parts be given. Then the following takes place
\begin{equation}
	n_F = 1_F\sum_{s=1}^n  \lambda_s(\vec{\beta}) \cdot (b_s)_F
\end{equation}
where 
\begin{equation}\label{eq:coeffmulti}
	\lambda_s(\vec{\beta}) = \lambda_K(b_s,b_{s+1}+...+b_k)\prod_{i=1}^{s-1} \lambda_M(b_i,b_{i+1}+...+b_k)
\end{equation}
or equivalent
\begin{equation}\label{eq:coeffmulti2}
	\lambda_s(\vec{\beta}) = \lambda_M(b_{s+1}+...+b_k,b_s)\prod_{i=1}^{s-1} \lambda_K(b_{i+1}+...+b_k,b_i).
\end{equation}
\end{corol}

\noindent {\it{\textbf{Proof.}}}
It is a straightforward algebraic induction exercise using  property (\ref{eq:TLambda}) of the sequence $\mathcal{T}_\lambda$.
The first form (\ref{eq:coeffmulti}) of the coefficients $\lambda_s(\vec{\beta})$ comes from the following 
$$
	\Big(b_1 + (n-b_1)\Big)_{\!F} \Rightarrow \Big(b_1 + b_2 + (n-b_1-b_2)\Big)_{\!F}
$$
while the second one (\ref{eq:coeffmulti2}) from
$$
	\Big((n-b_k) + b_k\Big)_{\!F} \Rightarrow \Big((n-b_k-b_{k-1}) + b_{k-1} + b_k\Big)_{\!F} \qquad \blacksquare
$$

\vspace{0.4cm}
If we take a vector $\langle 1,1,...,1\rangle$ of $n$ ones i.e. $b_s=1$ for any $s=1,2,...,n$; then we obtain alternative formula to compute elements of the sequence $F$.

\begin{corol}
Let $F\in\mathcal{T}_\lambda$ be given. Then $n$-th element of the sequence $F$ satisfies
\begin{equation}
	n_F = 1_F\cdot \sum_{s=1}^n \lambda_K(1,n-s)\prod_{i=1}^{s-1}\lambda_M(1,n-i)
\end{equation}
for any $n\in\mathbb{N}$.
\end{corol}

\begin{corol}
Let any sequence $F \in \mathcal{T}_\lambda$ be given. Then for any $n,k\in\mathbb{N}\cup\{0\}$ such that $n\geq k$, the $F$-nomial coefficients satisfy below recurrence identity
\begin{equation}\label{eq:rec}
	\fnomial{n}{k} = \lambda_K \fnomial{n-1}{k-1} + \lambda_M \fnomial{n-1}{k}
\end{equation}
\noindent where ${n \choose n}_F = {n \choose 0}_F = 1$.
\end{corol}

\noindent {\it{\textbf{Proof.}}} Take any $F\in\mathcal{T}_\lambda$ and $n\in\mathbb{N}\cup\{0\}$. Then from (\ref{eq:TLambda}) of $\mathcal{T}_\lambda$ and for any $m,k\in\mathbb{N}\cup\{0\}$ such that $m+k=n$ we have that $n$-th element of the sequence $F$ satisfies following recurrence
$$
	n_F = (k+m)_F = \lambda_K\cdot k_F + \lambda_M\cdot m_F
$$
\noindent Multiply both sides of above equation by $\frac{(n-1)_F!}{k_F!\cdot m_F!}$ to get

$$
\frac{n_F!}{k_F!\cdot m_F!} = \lambda_K\cdot \frac{(n-1)_F!}{(k-1)_F!\cdot m_F!} + \lambda_M\cdot \frac{(n-1)_F!}{k_F!\cdot(m-1)_F!}
$$

\noindent And from Definition \ref{def:fnomial} of $F$-nomial coefficients we have

$$
	\fnomial{n}{k} = \lambda_K {n-1 \choose k-1}_F + \lambda_M \fnomial{n-1}{k} \quad \blacksquare
$$

\vspace{0.4cm}
It turns out that the recurrence formula (\ref{eq:rec}) gives us a method to generating tilings of any layer $\layer{k}{n}$ designated by sequence $F\in\mathcal{T}_\lambda$.

\begin{theoremn}\label{th:1}
Let $F$ be a sequence of $\mathcal{T}_\lambda$ family. Then $F$ is cobweb tiling.
\end{theoremn}

\begin{figure}[ht]
\begin{center}
	\includegraphics[width=80mm]{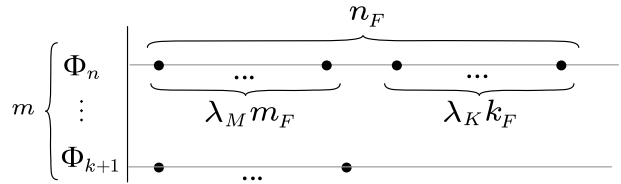}
	\caption{Picture of Theorem \ref{th:1} proof's idea. \label{fig:Step0}}
\end{center}
\end{figure}

\noindent {\it{\textbf{Proof.}}}
Suppose that we have a cobweb layer $\layer{k+1}{n}$ with $m$ levels designated by sequence $F$ from $\mathcal{T}_\lambda$ family and $m = n-k$. 
Consider $\Phi_n$ level with $n_F$ vertices. From (\ref{eq:TLambda}) we have that the number of vertices at this level is the sum of $\lambda_M \cdot m_F$ and $\lambda_K \cdot k_F$. Therefore we separate them by cutting into two disjoint subsets as illustrated by Figure \ref{fig:Step0} and cope at first $\lambda_M \cdot m_F$ vertices in Step 1. Then we shall cope the rest $\lambda_K \cdot k_F$ ones in Step~2.

\begin{figure}[ht]
\begin{center}
	\includegraphics[width=55mm]{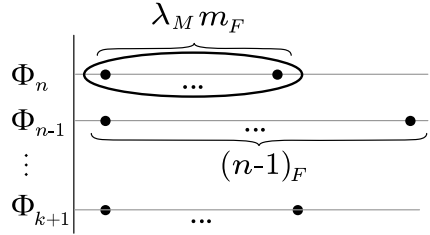}
	\caption{Picture of Theorem \ref{th:1} proof's Step 1. \label{fig:Step1}}
\end{center}
\end{figure}

\vspace{0.2cm}
{\it Step 1.}
Temporarily we have $\lambda_M \cdot m_F$ fixed vertices on $\Phi_n$ level to consider (Figure \ref{fig:Step1}). Let us cover them $\lambda_M$ times by $m$-th level of block $\sigma P_m$, which has exactly $m_F$ vertices. If $\lambda_M = 0$ we skip this step. What was left is the layer $\layer{k+1}{n-1}$ and we might eventually partition it with smaller disjoint blocks $\sigma P_{m-1}$ in the next induction step .

\begin{figure}[ht]
\begin{center}
	\includegraphics[width=110mm]{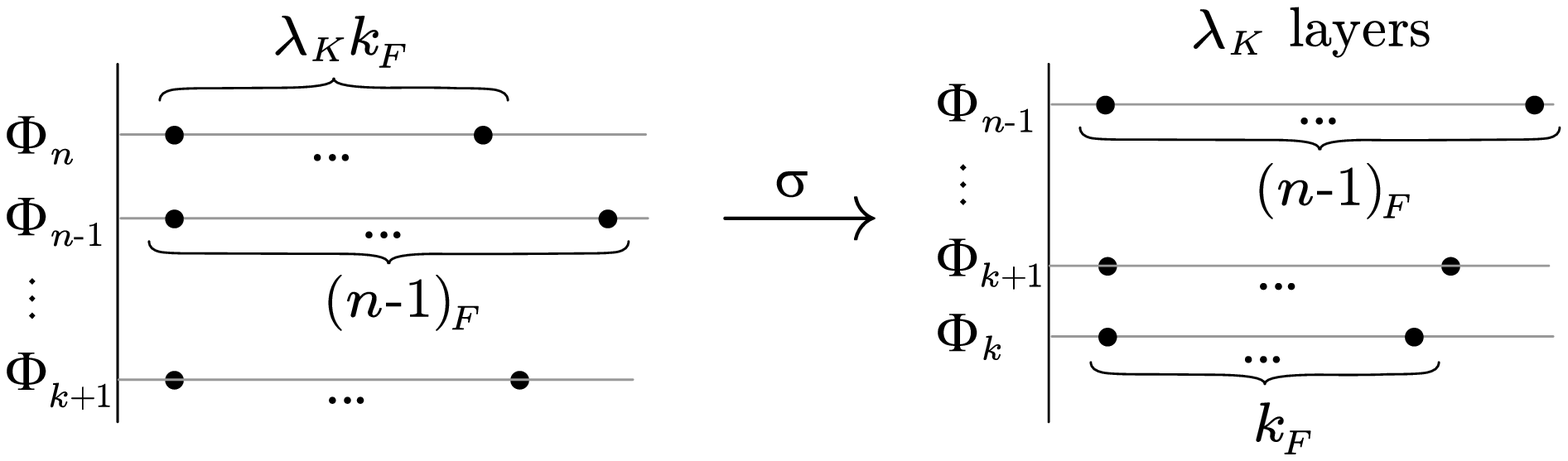}
	\caption{Picture of Theorem \ref{th:1} proof's Step 2. \label{fig:Step2}}
\end{center}
\end{figure}

{\it Step 2.}
Consider now the second complementary situation, where we have $\lambda_K \cdot k_F$ vertices on $\Phi_n$ level being fixed (Figure \ref{fig:Step2}). If $\lambda_K = 0$ we skip this step. Observe that if we move this level lower than $\Phi_{k+1}$ level, we obtain exactly $\lambda_K$ the same layers $\layer{k}{n-1}$ to be partitioned with disjoint blocks of the form $\sigma P_m$. This ``\emph{move}'' operation is just permutation $\sigma$ of levels' order.

\vspace{0.2cm}
{\it Recapitulation.}
The layer $\layer{k+1}{n}$ might be partitioned into $\sigma P_m$ blocks if $\layer{k+1}{n-1}$ might be partitioned into $\sigma P_{m-1}$ and $\layer{k}{n-1}$ into $\sigma P_m$ again. Continuing these steps by induction, we are left to prove that $\layer{k}{k}$ might be partitioned into $\sigma P_1$ blocks and $\layer{1}{m}$ into $\sigma P_m$ ones, what is trivial $\blacksquare$

\begin{observ}
	Let $F$ be a cobweb tiling sequence from the family $\mathcal{T}_\lambda$. Then the number $\Big\{ {n \atop k } \Big\}_F^1 $ of different tilings of layer $\layer{k}{n}$ where $n,k\in\mathbb{N}$, $n,k\geq 1$ is equal to:
\begin{equation} \label{eq:3}
	\bigg\{ {n \atop {k} } \bigg\}_F^1  =  \frac{n_F!}{(m_F!)^{\lambda_M} \cdot ((k-1)_F!)^{\lambda_K}}
	\cdot  \bigg( \bigg\{ {n-1 \atop k } \bigg\}_F^1 \bigg)^{\lambda_M}  	\cdot \bigg( \bigg\{ {n-1 \atop k-1 } \bigg\}_F^1 \bigg)^{\lambda_K}
\end{equation}
\noindent where $\Big\{ {n \atop n } \Big\}_F^1 = 1$ and $\Big\{ {n \atop 1 } \Big\}_F^1 = 1$.
\end{observ}

\noindent {\it{\textbf{Proof.}}}
According to steps of the proof of Theorem \ref{th:1} we might choose $m_F$ vertices $\lambda_M$ times at $n$-th level and next $(k-1)_F$ vertices $\lambda_K$ times out of $n_F$ ones in $\frac{n_F!}{(m_F!)^{\lambda_M} \cdot ((k-1)_F!)^{\lambda_K}}$ ways. Next recurrent steps of the proof of Theorem \ref{th:1} result in formula (\ref{eq:3}) via product rule of counting $\blacksquare$

\vspace{0.4cm}
Note that $\Big\{ {n \atop k } \Big\}_F^1$ is not the number of all different tilings of the layer  $\layer{k}{n}$ i.e. $\Big\{ {n \atop k } \Big\}_F^1 \leq \Big\{{n \atop k } \Big\}_F$  as computer experiments show \cite{md1}. There are much more other tilings with blocks $\sigma P_m$.


\section{Cobweb multi tiling}

In this section, more general case of the tiling problem is considered. For that to do we introduce the so-called multi $F$-nomial coefficients that counts blocks of multi-block partitions.

\begin{defn}\label{def:symbol}
	Let natural numbers' valued sequence $F\equiv\{n_F\}_{n\geq 0}$ and a composition $\langle b_1,b_2,...,b_k\rangle$ of the number $n$ be given. Then the multi $F$-nomial coefficient is identified with the symbol
\begin{equation}
	\fnomial{n}{b_1,b_2,...,b_k} = \frac{n_F!}{(b_1)_F!\cdot ... \cdot (b_k)_F!}
\end{equation}
while $n = b_1+b_2+...+b_k$.
\end{defn}

\begin{corol}
	Let $F$ be any $F$-cobweb admissible sequence. Then value of the multi $F$-nomial coefficient is natural number or zero i.e.
	\begin{equation}
	\fnomial{n}{b_1,b_2,...,b_k} \in \mathbb{N} \cup \{0\}
	\end{equation}
	for any $n,b_1,b_2,...,b_k\in\mathbb{N}$ such that $n=b_1+b_2+...+b_k$.
\end{corol}

For the sake of forthcoming combinatorial interpretation of multi $F$-nomial coefficients we introduce the following notation.

\begin{defn}
Let a cobweb layer $\layer{1}{n}$ of $n$ levels $\Phi_s$ and a composition $\langle b_1,b_2,...,b_k\rangle$ of number $n$ into $k$ non-zero parts be given. Then any cobweb layer $\langle\phi_1\to\phi_n\rangle$ of $n$ levels $\phi_s$ such that 
\begin{equation}
	\phi_s \subseteq \Phi_s, \qquad s=1,2,...,n;
\end{equation}
where the cardinality of $\phi_s$ is equal to $s$-th element of the vector $L$ given as follows
$$
	L = \sigma \cdot \langle 1, 2, ..., b_1, 1,2, ..., b_2, ..., 1,2, ..., b_k \rangle
$$
for any permutation $\sigma$ of a set $[n]$ is called cobweb multi-block of the form $\sigma P_{b_1,b_2,...,b_k}$. 
\end{defn}

\begin{figure}[ht]
\begin{center}
	\includegraphics[width=100mm]{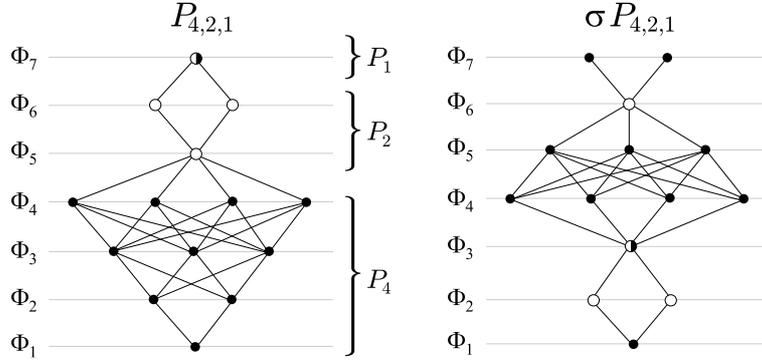}
	\caption{Examples of multi blocks $P_{4,2,1}$ and $\sigma P_{4,2,1}$. \label{fig:multiblock}}
\end{center}
\end{figure}

In the case of $\sigma = id$ we write for short $\sigma P_{b_1,b_2,...,b_k} = P_{b_1,b_2,...,b_k}$. Compare with Figure \ref{fig:multiblock}.

\begin{figure}[ht]
\begin{center}
	\includegraphics[width=110mm]{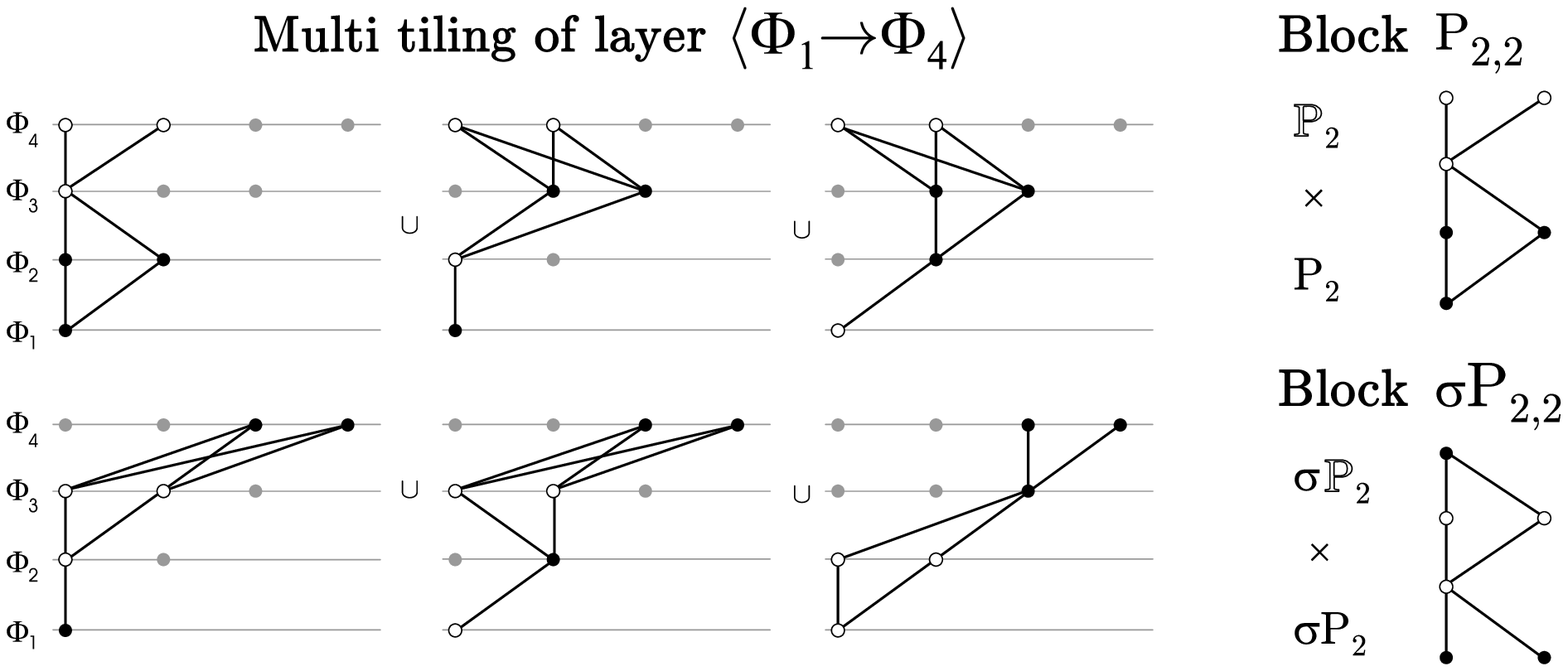}
	\caption{Sample multi tiling of layer $\layer{1}{4}$ from Example 2. \label{fig:sampletilingmulti}}
\end{center}
\end{figure}

\vspace{0.4cm}
\noindent \textbf{Example 1} \\
\noindent Take a sequence $F$ of next natural numbers i.e. $n_F = n$ and cobweb layer $\layer{1}{4}$ designated by $F$.
A sample multi tiling of the layer $\layer{1}{4}$ with the help of $\fnomial{4}{2,2} = 6$ disjoint multi blocks of the form $\sigma P_{2,2}$ is in Figure \ref{fig:sampletilingmulti}.

\begin{observ}
Let $\layer{1}{n}$ be a cobweb layer and $\langle b_1,...,b_k\rangle$ be a composition of the number $n$ into $k$ nonzero parts. 
Then the value of multi $F$-nomial coefficient $\fnomial{n}{b_1,b_2,...,b_k}$ is equal to the number of blocks that form the cobweb $\kappa$-partition, where $\kappa = |C_{max}(P_{b_1,...,b_k})|$. 
\end{observ}

\noindent {\it{\textbf{Proof.}}} The proof is natural extension of Observation 3 in \cite{akk1,akk2}.
The number of maximal paths in layer $\layer{1}{n}$ is equal to $n_F!$. However the number of maximal paths in any multi block $\sigma P_{b_1,b_2,...,b_k}$ is \break $\nobreak{(b_1)_F!\cdot(b_2)_F!\cdot...\cdot(b_k)_F!}$. Thus the number of such blocks is equal to 

$$
	\frac{n_F!}{(b_1)_F!\cdot(b_2)_F!\cdot...\cdot(b_k)_F!}
$$

\vspace{0.2cm}
\noindent where $n=b_1+b_2+...+b_k$ for any $n,k\in\mathbb{N}$ $\blacksquare$

\vspace{0.6cm}
\noindent Of course for $k=2$ we have

\begin{equation}
	{n \choose {b,n-b}}_F \equiv {n \choose b}_F = {n \choose {n-b}}_F
\end{equation}

\vspace{0.2cm}
\noindent \textbf{Note.} For any permutation $\sigma$ of the set $[k]$ the following holds

\begin{equation}
	{n \choose {b_1,b_2,...,b_k}}_F = {n \choose {b_{\sigma 1},b_{\sigma 2},...,b_{\sigma k}}}_F
\end{equation}

\vspace{0.2cm}
\noindent as is obvious from Definition \ref{def:symbol} of the multi F-nomial symbol. i.e.

$$
	\frac{n_F!}{(b_1)_F!\cdot(b_2)_F\cdot...\cdot(b_k)_F} = \frac{n_F!}{(b_{\sigma 1})_F!\cdot(b_{\sigma 2})_F\cdot...\cdot(b_{\sigma k})_F}
$$

\vspace{0.4cm}
\noindent
Let us observe also that for any natural $n,k$ and $b_1+...+b_m = n-k$ the following holds

\begin{equation} \label{eq:mult1}
	\fnomial{n}{k} \cdot \fnomial{n-k}{b_1,b_2,...,b_m} = \fnomial{n}{k,b_1,...,b_m}
\end{equation}

\begin{corol}
Let $F\in\mathcal{T}_\lambda$ and a composition $\vec{\beta}=\langle b_1,...,b_k\rangle$ of number $n$ into $k$ parts be given. Then the multi $F$-nomial coefficients satisfy the following recurrence relation
\begin{equation}
	\fnomial{n}{b_1,b_2,...,b_k} = \sum_{s=1}^k \lambda_s(\vec{\beta}) \cdot
	\fnomial{n-1}{b_1,...,b_{s-1},b_s-1,b_{s+1},...,b_k}
\end{equation}
for coefficients $\lambda_s(\vec{\beta})$ from (\ref{eq:coeffmulti}) and for any  $n= b_1 +...+b_k$ and $\fnomial{n}{n, 0, ... ,0} = 1$.
\end{corol}

\noindent {\it{\textbf{Proof.}}}
Take any $F\in\mathcal{T}_\lambda$ and a composition $\vec{\beta}=\langle b_1,...,b_k\rangle$ of the number $n$. Then from Corollary \ref{cor:lambdamult} we have that for certain coefficients $\lambda_s(\vec{\beta})$ any $n$-th element of the sequence $F$ satisfies
$$
	n_F = \sum_{s=1}^k \lambda_s(\vec{\beta}) \cdot (b_s)_F
$$
\noindent If we multiply both sides by $\frac{(n-1)_F!}{(b_1)_F!\cdot ... \cdot(b_k)_F!}$ then we obtain
$$
	\fnomial{n}{b_1,...,b_k} = \sum_{s=1}^k \lambda_s(\vec{\beta})
	\frac{(n-1)_F!}{(b_1)_F!\cdot...\cdot(b_{s-1})_F!(b_s-1)_F!(b_{s+1})_F!\cdot ...\cdot(b_k)_F!}
$$
Hence the thesis $\blacksquare$

\begin{theoremn}\label{th:multi}
	Let any sequence $F \in \mathcal{T}_\lambda$ be given. Then the sequence $F$ is cobweb multi tiling i.e. any layer $\layer{1}{n}$ might be partitioned into multi-blocks of the form $\sigma P_{b_1,b_2,...,b_k}$ such that $b_1+...+b_k=n$.
\end{theoremn}

\noindent {\it{\textbf{Proof.}}}
Take any cobweb layer $\layer{1}{n}$ designated by sequence $F\in\mathcal{T}_\lambda$ and a number $k\in\mathbb{N}$. We need to partition the layer into disjoint multi blocks of the form $\sigma P_{b_1,b_2,...,b_k}$.

\begin{figure}[ht]
\begin{center}
	\includegraphics[width=80mm]{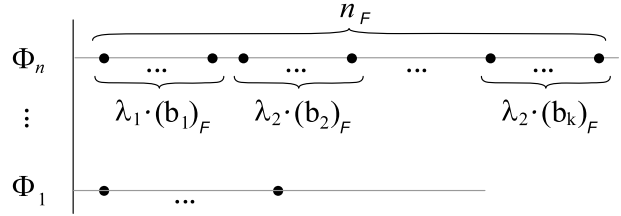}
	\caption{Idea's picture of Theorem \ref{th:multi}. \label{fig:proof}}
\end{center}
\end{figure}

\vspace{0.2cm}
\noindent Consider level $\Phi_n$ with $n_F$ vertices. 
From Corollary \ref{cor:lambdamult} we have that the number of vertices at this level is the following sum
$$
	n_F = \sum_{s=1}^{k}{\lambda_s(\vec{\beta}) \cdot (b_s)_F}
$$
for certain coefficients $\lambda_s(\vec{\beta})$ where $1\leq s \leq k$ and $\vec{\beta} = \langle b_1,b_2,...,b_k\rangle$.

\vspace{0.2cm}
\noindent Therefore let us separate these $n_F$ vertices by cutting into $k$ disjoint subsets as illustrated by Fig. \ref{fig:proof} and cope at first $\lambda_1\cdot(b_1)_F$ vertices in Step 1, then $\lambda_2\cdot(b_2)_F$ ones in Step 2 and so on up to the last $\lambda_k\cdot(b_k)_F$  vertices to consider in the last $k$-th step. If any $\lambda_i = 0$ we skip $i$-th step.

\vspace{0.2cm}
{\it Step 1.}
Temporarily we have $\lambda_1\cdot(b_1)_F$ fixed vertices at level $\Phi_n$ to consider. Let us cover them $\lambda_1$ times by $(b_1)$-th level of block $P_{b_1,b_2,...,b_k}$, which has exactly $(b_1)_F$ vertices. What was left is the layer $\layer{1}{n-1}$ and we might partition it with smaller disjoint blocks $\sigma P_{b_1-1,b_2,...,b_k}$ in the next induction step. 

\vspace{0.2cm}
{Note.} In the next induction steps we use smaller blocks $\sigma P$ without levels which we have been already  used in previous steps (disjoint of blocks condition).

\vspace{0.2cm}
{\it Step 2.}
Consider now the second situation, where we have $\lambda_2\cdot(b_2)_F$ vertices at level $\Phi_n$ being fixed. We cover them $\lambda_2$ times by $(b_1+b_2)$-th level of block $P_{b_1,b_2,...,b_k}$, which has $(b_2)_F$ vertices. Then we obtain smaller layer $\layer{1}{n-1}$ to be partitioned with blocks $\sigma P_{b_1,b_2-1,b_3,...,b_k}$.

\vspace{0.2cm}
\noindent And so on up to ...

\vspace{0.2cm}
{\it Step $k$.}
Analogously to previous steps, we cover the last $\lambda_{b_s}$ vertices by the last $(b_1+b_2+...+b_k)=n$-th level of block $P_{b_1,b_2,...,b_k}$, obtaining smaller layer $\layer{1}{n-1}$ to be partitioned with blocks $\sigma P_{b_1,...,b_{k-1},b_k-1}$.

\vspace{0.2cm}
{\it Conclusion.}

\noindent The layer $\layer{1}{n}$ might be partitioned into blocks $\sigma P_{b_1,b_2,...,b_k}$ if  $\layer{1}{n-1}$ might be partitioned into $\sigma P_{b_1-1,b_2,...,b_k}$ and $\layer{1}{n-1}$ into $\sigma P_{b_1,b_2-1,b_3,...,b_k}$ again and so on up to the layer $\layer{1}{n-1}$ which might be partitioned into $\sigma P_{b_1,...,b_{k-1},b_k-1}$. Continuing these steps by induction, we are left to prove that $\layer{1}{k}$  might be partitioned into blocks $\sigma P_{1,1,...,1}$ or $\layer{1}{1}$ by $\sigma P_{1,0,...,0}$ ones, which is trivial. $\blacksquare$

\section{Family $\mathcal{T}_\lambda(\alpha,\beta)$ of cobweb tiling sequences}

In this section a specific family of cobweb tiling sequences $F\in\mathcal{T}_\lambda$ is presented as an exemplification of a might be source method.
We assume that coefficients $\lambda_K$ and $\lambda_M$ of $F\in\mathcal{T}_\lambda$ take a form
\begin{equation}
	\lambda_M(k,m) = \alpha^k \qquad 
	\lambda_K(k,m) = \beta^m 
\end{equation}
while $\alpha,\beta\in\mathbb{N}$.

\begin{notation}
Let $\mathcal{T}_\lambda(\alpha,\beta)$ denotes a family of natural numbers' valued sequences $F\equiv\{n_F\}_{n\geq 0}$ constituted by $n$-th coefficients of the generating function $\mathcal{F}(x)$ expansion i.e. $n_F = [x^n]\mathcal{F}(x)$, where
\begin{equation}\label{eq:form}
	\mathcal{F}(x) = 1_F\cdot \frac{x}{(1-\alpha x)(1-\beta x)}
\end{equation}
for certain $\alpha,\beta\in \mathbb{N}\cup\{0\}$ and $1_F \in \mathbb{N}$.
\end{notation}

\begin{enumerate}
	\item If ($\alpha = \beta$), then $\mathcal{F}(x) = 1_F\cdot \frac{x}{1-\alpha x} + \alpha x \mathcal{F}(x)$ which leads to
	\begin{equation}\label{eq:aa}
		n_F = 1_F\cdot n \cdot \alpha^{n-1} \qquad n\geq 1
	\end{equation}
	
	\item  If ($\alpha \neq \beta$), then $\mathcal{F}(x) = \frac{1_F}{\alpha-\beta}\left( \frac{1}{1-\alpha x} - \frac{1}{1 - \beta x} \right)$ gives us
	
	\begin{equation}\label{eq:ab}
		n_F = \frac{1_F}{\alpha - \beta}\left( \alpha^n - \beta^n \right) \qquad n\geq 1
	\end{equation}
\end{enumerate}

\begin{proposition}\label{prop:sum}
Let $F\in\mathcal{T}_\lambda(\alpha,\beta)$ and composition $\vec{b} = \langle b_1,b_2,...,b_k\rangle$ of the number $n$ into $k$ non-zero parts be given. Then any $n$-th element of the sequence $F$ satisfies the following recurrence identity

\begin{equation}\label{eq:sum}
	n_F = \left( \sum_{s=1}^{k} b_s \right)_{\!\!F} = \sum_{s=1}^k \lambda_s(\vec{b}) \cdot (b_s)_F
\end{equation}
where
$$
	\lambda_s(\vec{b}) = \alpha^{b_{s+1} + ... + b_{k}}\cdot \beta^{b_1+...+b_{s-1}}
$$
for any $n=b_1+...+b_k$.
\end{proposition}

\noindent {\it{\textbf{Proof.}}} Take any composition $\vec{b} = \langle b_1,b_2,...,b_k\rangle$ of the number $n\in\mathbb{N}$ into $k$ nonzero parts i.e. $b_1+b_2+...+b_k=n$.
\begin{enumerate}
\item If ($\alpha = \beta$) then from (\ref{eq:aa})


$ \left( \sum_{s=1}^k b_s \right)_{\!\!F} = 1_F \left( \sum_{s=1}^k b_s \right) \cdot \alpha^{n-1} = \sum_{s=1}^k 1_F b_s \alpha^{b_s-1} \alpha^{n-b_s} = $

$ = \sum_{s=1}^k (b_s)_F \alpha^{n-b_s}$

\item If ($\alpha \neq \beta$) then from (\ref{eq:ab})


$\left( \sum_{s=1}^k b_s \right)_{\!\!F} = \frac{1_F}{\alpha-\beta}\alpha^{b_1+\sum_{s=2}^k b_s} - \frac{1_F}{\alpha-\beta}\beta^{b_k + \sum_{s=1}^{k-1} b_s } = A + B $

\noindent Next, denote $S_{\pm}(m)$ for $1<m<k$ such that $S_{+}(m) + S_{-}(m) = 0$ as follows
$ S_{\pm}(m) = \pm\frac{1_F}{\alpha-\beta} \alpha^{\sum_{s=m+1}^k b_s} \cdot \beta^{\sum_{s=1}^m b_s}$.
\noindent Then observe that if we add to the $A+B$ the sum of $S_{\pm}(m)$ where $1<m<k$ i.e.

$ A + B = A + B + \sum_{1<j<k} S_{+}(j) + S_{-}(j) $
\noindent then we obtain

$
\left\{
\begin{array}{l}
A + S_{-}(1) = (b_1)_F \cdot \alpha^{\sum_{s=2}^n b_s} \beta^{0} \\
S_{+}(1) + S_{-}(2) = (b_2)_F \cdot \alpha^{\sum_{3=2}^n b_s} \cdot \beta^{b_1} \\
... \\
S_{+}(k-1) + B = (b_k)_F \cdot \alpha^{0} \cdot \beta^{\sum_{s=1}^{k-1} b_s} \\
\end{array}
\right.
$

\noindent And finally \\
$ \left( \sum_{s=1}^k b_s \right)_{\!\!F} = A + B = \sum_{s=1}^k (b_s)_F\cdot \alpha^{b_{s+1}+...+b_k}\beta^{b_1+...+b_{s-1}}$ $\blacksquare$
\end{enumerate}

\vspace{0.2cm}
\noindent \textbf{Note.} If $k = 2$ then for any $m,b\in\mathbb{N}\cup\{0\}$ we have
\begin{equation}\label{eq:two}
	(m+b)_F = \lambda_M m_F + \lambda_b b_F = \alpha^b m_F + \beta^m b_F
\end{equation}

Let us compare above with condition (\ref{eq:TLambda}) for sequences that are cobweb tiling 
from family $\mathcal{T}_\lambda$ and let us sum up this with the following corollary.

\begin{corol}
Let family of sequences $\mathcal{T}_\lambda(\alpha,\beta)$ and family $\mathcal{T}_\lambda$ of cobweb tiling sequences be given. Then the following takes place
\begin{equation}
	\mathcal{T}_\lambda(\alpha,\beta) \subset \mathcal{T}_\lambda
\end{equation}
thus any sequence $F\in\mathcal{T}_\lambda(\alpha,\beta)$ is cobweb tiling. 
\end{corol}

\noindent {\it{\textbf{Proof.}}}
We only need to show that $\mathcal{T}_\lambda(\alpha,\beta) \neq \mathcal{T}_\lambda$. 
As an example we show that the sequence $F$ of Fibonacci numbers is cobweb tiling of the 
form $\mathcal{T}_\lambda$ but does not belong to the family $\mathcal{T}_\lambda(\alpha,\beta)$.
Ones show that $n$-th element of the Fibonacci numbers satisfies
\begin{equation}
	n_F = \frac{1}{\alpha - \beta}\left( \alpha^n - \beta^n \right)
\end{equation}
but $\alpha = \frac{1 + \sqrt{5}}{2}$ and $\beta = \frac{1-\sqrt{5}}{2}$ are not natural numbers 
- compare with (\ref{eq:form}).
However its elements satisfy another equivalent relation for any $m,k\in\mathbb{N}\cup\{0\}$
\begin{equation}
	(k+m)_F = (m-1)_F\cdot k_F + (k+1)_F\cdot m_F
\end{equation}
\vspace{0.2cm}
\noindent Therefore $F\in\mathcal{T}_\lambda$ and  $F\notin\mathcal{T}_\lambda(\alpha,\beta)$. Hence the thesis $\blacksquare$

\begin{corol}
Let $F\in\mathcal{T}_\lambda$ be given. Then for any $n,k\in\mathbb{N}\cup\{0\}$ the following holds
\begin{equation}
	(k\cdot n)_F = \bigg( \underbrace{n + n + ... + n}_k \bigg)_{\!\!F} = n_F \cdot \sum_{s=1}^k \alpha^{(k-s)n} \beta^{(s-1)n}
\end{equation}
\end{corol}

\vspace{0.4cm}
From Proposition \ref{prop:sum} we obtain an another explicit formula for $n$-th element of the sequence $F\in\mathcal{T}_\lambda$ i.e.

\begin{equation}
	n_F = (n \cdot 1)_F = 1_F\sum_{s=1}^n \alpha^{(n-s)} \beta^{(s-1)}.
\end{equation}

\section{Examples of cobweb tiling sequences} \label{sect:examples}

In this section we are going to show a few examples of cobweb-tiling sequences.
Throughout this part we shall consequently use the condition convention:  $n = k + m$.

\subsection{Examples of $\mathcal{T}_\lambda(\alpha,\beta)$ family }

\begin{enumerate}
\item \textbf{Natural numbers} \\
	Putting $\alpha=\beta= 1$ gives us a sequence $n_F = 1_F\cdot n$ with the recurrence $(k+m)_F = k_F + m_F$. If $1_F=1$ then we obtain Natural numbers with Binomial coefficients' recurrence:
$$
	{n \choose k} \equiv \fnomial{n}{k} = \fnomial{n-1}{k-1} + \fnomial{n-1}{k}
$$

\item \textbf{Powers' sequence} \\
	If $\alpha = 0, \beta = 1_F = q$ then $n_F = q^n$ and $(k+m)_F = q^m\cdot k_F$ with its $F$-nomial coefficients' recurrence
$$
	\fnomial{n}{k} = q^m \fnomial{n-1}{k-1} = q^k \fnomial{n-1}{m-1}
$$

\item \textbf{Gaussian numbers} \\
	If $\alpha = 1, \beta = q$ then $n_F = \frac{1_F}{1-q}\left( 1 - q^n \right)$ and  $(k+m)_F = k_F + q^k m_F$ with the recurrence for Gaussian coefficients
$$
	\fnomialF{n}{k}{q} \equiv \fnomial{n}{k} = \fnomial{n-1}{k-1} + q^k \fnomial{n-1}{k}
$$

\item \textbf{Modified Gaussian integers} \label{ex:modgaus}\\
	For $\alpha = \beta = q\in\mathbb{N}$ we have $n_F = 1_F\cdot n \cdot q^{n-1}$ and  $(k+m)_F = q^m k_F + q^k m_F$ with the recurrence
$$
	\fnomial{n}{k} = q^m \fnomial{n-1}{k-1} + q^k \fnomial{n-1}{k}
$$

\end{enumerate}

\subsection{Fibonacci numbers}

In the following, we prove that sequence of Fibonacci numbers is tiling sequence i.e. any cobweb layer $\layer{k}{n}$ might be partitioned into blocks of the form $\sigma P_m$.

\begin{defn}
Let $F(p)$ be a natural numbers' valued sequence such that for any $k,m\in\mathbb{N}\cup\{0\}$ its elements satisfy the following relation
\begin{equation}\label{eq:fib}
	(k+m)_F = (m-1)_F\cdot k_F + (k+1)_F\cdot m_F
\end{equation}
while $1_F = 1$ and $2_F = p$.
\end{defn}

\vspace{0.2cm}
\noindent From Theorem \ref{th:1} and condition (\ref{eq:TLambda}) on the sequence $\mathcal{T}_\lambda$, we have that $F(p)$ is cobweb tiling. 
Moreover, it is easy to see, that explicit formula for $n$-th element of $F(p)$ is 
\begin{equation}
	n_F = \frac{1}{\sqrt{2_F^2 + 4}}\left( \phi_1^n - \phi_2^n \right)
\end{equation}
where $\phi_{1,2} = \frac{2_F\pm\sqrt{2_F^2+4}}{2}$ and $1_F = 1$ while $2_F = p$.

\vspace{0.4cm}
\noindent \textbf{Examples of $F(p) = \{n_F\}_{n\geq 0}$}
\begin{itemize}
\item $F(1) \equiv (0, 1, 1, 2, 3, 5, 8, 13, 21, 34, 55, 89, 144, ...) \equiv$ Fibonacci numbers
\item $F(2) \equiv (0, 1, 2, 5, 12, 29, 70, 169, 408, 985, 2378, 5741, 13860,...)$
\item $F(3) \equiv (0, 1, 3, 10, 33, 109, 360, 1189, 3927, 12970, 42837,...)$
\item $F(4) \equiv ( 0, 1, 4, 17, 72, 305, 1292, 5473, 23184, 98209, 416020,...)$

\end{itemize}

\begin{corol}
The sequence of Fibonacci numbers is cobweb tiling.
\end{corol}

\noindent {\it{\textbf{Proof.}}}
If we put $1_F=2_F=1$ in (\ref{eq:fib}) then we obtain \textbf{Fibonacci numbers} and well-known recurrence relation for Fibonomial coefficients \cite{eks}

\begin{equation}
	\fnomial{n}{k} = (m-1)_F\fnomial{n-1}{k-1} + (k+1)_F\fnomial{n-1}{k} \blacksquare
\end{equation}

\begin{observ}
Let $F$ be a sequence of the form $F(p)$. Take any composition $\langle b_1,b_2,...,b_k \rangle$ of a number $n$ into $k$ nonzero parts. Then $n$-th element of $F$ satisfies
\begin{equation}
	n_F = \sum_{s=1}^k (b_s)_F\cdot\prod_{i=1}^{s-1}(b_i+1)_F\cdot(b_{s+1}+...+b_k-1)_F
\end{equation}
while $n,k\in\Nat$.
\end{observ}

\noindent {\it{\textbf{Proof.}}}
It is a straightforward algebraic exercise using an idea from the proof of Corollary \ref{cor:lambdamult}.
If we use the substitutions $m = a + b$ in the formula (\ref{eq:fib}) then we obtain the case of $3$ terms
$$
	(k + m)_F = (k + a + b)_F = \lambda_K k_F + \lambda_a a_F + \lambda_b b_F
$$
where $\lambda_K\! =\! (a+b-1)_F$, $\lambda_a\! =\! (k + 1)_F\cdot(b - 1)_F$ 
and $\lambda_b\! =\! (k+1)_F\cdot(a+1)_F$. And so on by induction $\blacksquare$

\section{Cobweb tiling problem as a particular case of clique problem}

Recall that the clique problem is the problem of determining whether a graph contains a clique of at least a given size $d$. In this section, we show that the cobweb tiling problem might be considered as the clique problem in specific graph.
Namely reformulation of the $F$-cobweb i.e. $F$-boxes tiling problem into a clique problem of a graph specially invented for that purpose - is proposed.

\vspace{0.4cm}
Suppose that we have a cobweb layer $\layer{k}{n}$ designated by any sequence $F$. Let $B\left(\layer{k}{n}\right)$ denotes a family of all blocks of the form $\sigma P_m$, where $m=n-k+1$ of that layer $\layer{k}{n}$ and assume that $b_{k,n}$ is a 
cardinality of that family i.e. $b_{k,n} = |B\left(\layer{k}{n}\right)|$.

\begin{observ}
	The number $b_{k,n}$ is given by the following formula
$$
	b_{k,n} = \sum_{\sigma\in S_m}{\prod_{s=1}^{m}{{ 
	{(k+s-1)_F}
	\choose
	{(\sigma\cdot s)_F} 
	}}}
$$
where $m=n-k+1$ and $S_m$ is a set of permutations $\sigma$ of the set $\break\{k_F,(k+1)_F,...,n_F\}$.
\end{observ}

\noindent {\it{\textbf{Proof.}}}
\noindent Suppose that we have the layer $\layer{k}{n}$. Take any permutation $\sigma\in S_m$ of $m$ 
levels of the block $\sigma P_m$. Let $s\in[m]$; for such order of levels, cope $(\sigma\cdot s)_F$ 
vertices by $s$-th element of the block $\sigma P_m$ from all of vertices i.e.
$(k+s-1)_F$ of the $(k+s)$-th level in the layer $\layer{k}{n}$. To the end sum the above after 
all of permutation $\sigma$ $\blacksquare$

\vspace{0.4cm}
Let us define now a simple not directed graph $G(\layer{k}{n})=(V,E)$ such that set of vertices is $V\equiv B\left(\layer{k}{n}\right)$ i.e. for any cobweb block $\beta$ we have that
$$
	\beta\in B\left(\layer{k}{n}\right) \Leftrightarrow v_\beta\in V 
$$
while set of edges $E$ is defined as follows
$$
	\{ v_\alpha, v_\beta \} \in E \Leftrightarrow C_{max}(\alpha) \cap C_{max}(\beta) = \emptyset
$$
\noindent for any two cobweb blocks $\alpha,\beta \in B\left(\layer{k}{n}\right)$ where $C_{max}(\gamma)$ is a set of maximal paths of block $\gamma$.

\begin{corol}
	Cobweb tiling problem of layer $\layer{k}{n}$ is the clique of size $d$  in graph $G(\layer{k}{n})$  problem , where $d = m_F!$.
\end{corol}

\noindent {\it{\textbf{Proof.}}}
\noindent Suppose that we have a cobweb layer $\layer{k}{n}$ and consider the family $B\left(\layer{k}{n}\right)$ of all blocks of the form $\sigma P_m$ of layer $\layer{k}{n}$, where $m=n-k+1$. 

Assume that a cobweb tiling of layer $\layer{k}{n}$ contains $d$ pairwise disjoint blocks of the form $\sigma P_m$, where $m=n-k+1$. From combinatorial interpretation of $F$-nomial coefficients we have that $d = \fnomial{n}{m}$. Thus if the family $B\left(\layer{k}{n}\right)$ contains $d$ blocks that are pairwise disjoint then the layer has tiling $\pi$. In the other words, if a graph $G$ has $d$ vertices that are pairwise incidence then of course has a clique $\chi$ of size $d$. Moreover this clique $\chi$ of graph $G$ corresponds to the cobweb tiling $\pi$ of layer $\layer{k}{n}$ and vice versa i.e. $\pi \Leftrightarrow \chi$ $\blacksquare$

\begin{corol}
If a graph $G(\layer{k}{n})$ has a clique $\chi$ of size $d=m_F!$ then $\chi$ is maximal clique of the graph.
\end{corol}

\begin{corol}
The number of all cobweb tilings of layer $\layer{k}{n}$ is equal to the number of all maximal cliques in graph $G(\layer{k}{n})$.
\end{corol}

\section{Map of cobweb sequences}

\begin{figure}[ht]
\begin{center}
	\includegraphics[width=83mm]{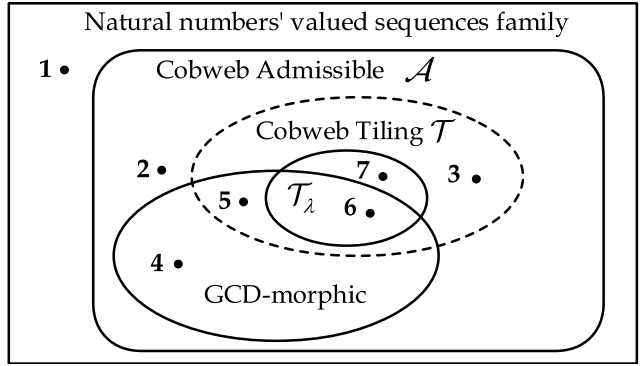}
	\caption{Venn type map of various cobweb sequences families. \label{fig:CobwebMap}}
\end{center}
\end{figure}

Here down in Figure \ref{fig:CobwebMap} we present a Venn type diagram map of cobweb sequences. Note that the boundary of the whole family of Cobweb Tilling sequences is still not known (open problem).

\vspace{0.2cm}
\textbf{Cobweb Admissible} sequences family $\mathcal{A}$ is defined in \cite{md2}, 
\textbf{GCD-morphic} sequences family in \cite{md1}. 
Subfamily $\mathcal{T}_\lambda$ of \textbf{cobweb tiling} sequences $\mathcal{T}$ is introduced in this note.

\begin{enumerate}
\item $A = (1,3,5,7,9,...) $;
\item $B = (1,2,2,2,1,4,1,2,...) = B_{2,2} \cdot B_{2,3}$;
\item $C = (1,2,2,1,2,2,1,...)$;
\item $E = (1,2,3,2,1,6,1,...) = B_{2,2} \cdot B_{3,3}$;
\item $F = (1,2,1,2,1,2,...) = B_{2,2} $;
\item Natural numbers, Fibonacci numbers;
\item $G = 1, 4, 12, 32, 80, 192, 448, 1024, ...  $ (Example \ref{ex:modgaus} in Section \ref{sect:examples});
\end{enumerate}

\noindent Sequences $B_{c,M}$ and $A_{c,t}$ are defined in \cite{md1}.

\vspace{0.4cm}
\noindent \textbf{Additional information}

\vspace{0.2cm}
In \cite{mdweb} we deliver some computer applications for generating tilings of any layer $\layer{k}{n}$ based on an algorithm from the proof of Theorem \ref{th:1}. There one may find also a visualization application for drawing all multi blocks of the form $\sigma P_{k,n-k}$ of a layer $\layer{1}{n}$.

\vspace{0.4cm}
\noindent \textbf{Acknowledgments}

\vspace{0.2cm}
I would like to thank my Professor A. Krzysztof Kwa\'sniewski for his  comments and effective improvements of this note.



\begin{thebibliography}{99}

\bibitem{akk1}  A. Krzysztof Kwa\'sniewski, {\it On cobweb posets and their combinatorially admissible sequences},
Adv. Studies Contemp. Math. Vol. 18 (1), 2009, ArXiv:math/0512578, 19 Jan 2009  

\bibitem{akk2}  A. Krzysztof Kwa\'sniewski, {\it Cobweb posets as noncommutative prefabs},
Adv. Stud. Contemp. Math.  vol.14 (1) 2007. pp. 37-47; ArXiv: math.CO/0503286, 25 Sep 2005

\bibitem{akk3} A. Krzysztof Kwa\'sniewski, {\it Fibonomial cumulative connection constants} 
ArXiv:math/0406006, 20 Feb 2009, upgrade of Bulletin of the ICA vol. 44 (2005) 81-92.

\bibitem{akk4} A. K. Kwa\'sniewski , {\it Graded posets inverse zeta matrix  formula}
ArXiv:0903.2575, 14 Mar 2009 (130-th  Birthday of Albert Einstein)


\bibitem{akkmd1} A. Krzysztof Kwa\'sniewski, M. Dziemia\'nczuk, {\it Cobweb posets - Recent Results},
ISRAMA 2007, December 1-17 2007 Kolkata, INDIA, Adv. Stud. Contemp. Math. volume 16 (2), 2008 (April) pp. 197-218
ArXiv:0801.398, 23 Mar 2009

\bibitem{akkmd2}  A. K. Kwa\'sniewski, M. Dziemia\'nczuk  {\it On cobweb posets' most relevant codings}, Preprint ArXiv:0804.1728, 27 Feb 2009


\bibitem{eks} E. Krot, {\it An Introduction to Finite Fibonomial Calculus}, 
CEJM 2(5) (2005) 754-766.

\bibitem{md1} M. Dziemia\'nczuk, {\it On Cobweb posets tiling problem},
Adv. Stud. Contemp. Math. volume 16 (2), 2008 (April) pp. 219-233
ArXiv:0709.4263, 4 Oct 2007

\bibitem{md2} M. Dziemia\'nczuk, \emph{On multi F-nomial coefficients and Inversion formula for F-nomial coefficients}, 
ArXiv:0806.3626, 23 Dec 2008

\bibitem{md3} M. Dziemia\'nczuk, \emph{On Cobweb Admissible Sequences - The Production Theorem},
Proceedings of The 2008 International Conference on Foundations of Computer Science (FCS'08),
July 14-17, 2008, Las Vegas, USA pp.163-165, ArXiv:0801.4699, 30 Jan 2008

\bibitem{md4} M. Dziemia\'nczuk, W.Bajguz, \emph{On GCD-morphic sequences},
ArXiv:0802.1303, 10 Feb 2008

\bibitem{mdweb} M.Dziemia\'nczuk, \emph{Cobweb Posets Website},\\
http://www.faces-of-nature.art.pl/cobwebposets.html

\end{thebibliography}
\end{document}